\documentclass{article}
\usepackage{graphicx}
\usepackage{tikz}
\usepackage{geometry}
\usepackage{setspace}
\usepackage{indentfirst}
\usepackage{xcolor}
\usepackage{hyperref}
\usepackage{amsmath}
\usepackage{amsfonts}
\usepackage{float}
\usepackage{amssymb,amsthm,graphicx}

\hypersetup{
    colorlinks=true,
    linkcolor=blue,
    citecolor=blue,
    filecolor=blue,
    urlcolor=blue,
}

\newtheorem{theorem}{Theorem}[section]
\newtheorem{lemma}[theorem]{Lemma}

\title{\textbf{The Stability of Persistence Diagrams Under Non-Uniform Scaling} \\[0.5em]}
\author{
\textbf{Vu Anh Le\textsuperscript{1*}, Mehmet Dik\textsuperscript{1,2**}}  \\
\textsuperscript{1} Department of Mathematics and Computer Science, Beloit College \\
\textsuperscript{2} Department of Mathematics, Computer Science \& Physics, Rockford University \\ [1em]
\textit{Contact} \\
* \href{mailto:csplevuanh@gmail.com}{csplevuanh@gmail.com} \\
** \href{mailto:mdik@rockford.edu}{mdik@rockford.edu}
}
\date{\today}

\begin{document}

\maketitle

\begin{abstract}
We investigate the stability of persistence diagrams \( D \) under non-uniform scaling transformations \( S \) in \( \mathbb{R}^n \). Given a finite metric space \( X \subset \mathbb{R}^n \) with Euclidean distance \( d_X \), and scaling factors \( s_1, s_2, \ldots, s_n > 0 \) applied to each coordinate, we derive explicit bounds on the bottleneck distance \( d_B(D, D_S) \) between the persistence diagrams of \( X \) and its scaled version \( S(X) \). Specifically, we show that 
\[
d_B(D, D_S) \leq \frac{1}{2} (s_{\max} - s_{\min}) \cdot \operatorname{diam}(X),
\]
where \( s_{\min} \) and \( s_{\max} \) are the smallest and largest scaling factors, respectively, and \( \operatorname{diam}(X) \) is the diameter of \( X \). We extend this analysis to higher-dimensional homological features, alternative metrics such as the Wasserstein distance, and iterative or probabilistic scaling scenarios. Our results provide a framework for quantifying the effects of non-uniform scaling on persistence diagrams. 
\end{abstract}

\begin{center}
\tableofcontents
\end{center}

\vspace{0.2in}

\section{Introduction}

Persistent homology is a central tool in topological data analysis (TDA), providing a robust framework for summarizing topological features of datasets via persistence diagrams \( D \) \cite{ZC05, EH08}. A key property of persistent homology is its stability under small perturbations, ensuring that minor changes in the input data lead to controlled variations in \( D \) \cite{CSEH07, CDS10}.

Non-uniform scaling, where each coordinate is scaled by a distinct factor, arises frequently in applications such as image processing \cite{JSZK15} and high-dimensional feature normalization \cite{IS15}. Unlike uniform scaling or isometric transformations, non-uniform scaling can introduce anisotropic distortions, significantly affecting inter-point distances and topological features. While prior work has addressed stability under uniform scaling and small perturbations \cite{CSEH07, BL15}, the effects of non-uniform scaling remain underexplored.

This paper addresses this gap by analyzing the stability of persistence diagrams under non-uniform scaling transformations \( S \), defined as
\[
S(x_1, x_2, \ldots, x_n) = (s_1 x_1, s_2 x_2, \ldots, s_n x_n),
\]
where \( s_1, s_2, \ldots, s_n > 0 \). Our contributions include deriving explicit bounds on the bottleneck distance \( d_B(D, D_S) \), such as:
\[
d_B(D, D_S) \leq \frac{1}{2} (s_{\max} - s_{\min}) \cdot \operatorname{diam}(X),
\]
where \( s_{\min} \) and \( s_{\max} \) are the smallest and largest scaling factors.
We then extend stability results to higher homology dimensions, alternative filtrations, and metrics like the Wasserstein distance. We also investigate the iterative and probabilistic scaling scenarios to quantify cumulative and expected impacts.

Our findings provide theoretical guarantees and practical insights into the behavior of persistence diagrams under non-uniform scaling, with applications to noisy, high-dimensional, and anisotropically scaled data.

\vspace{0.3in}

\section{Theoretical Framework}
\subsection{Persistence Diagrams}
Given \( X \subset \mathbb{R}^n \), we construct a filtration \( \{K_\epsilon\}_{\epsilon \geq 0} \) using methods such as the Vietoris--Rips or Čech complexes. Persistent homology studies the birth and death of homological features as \( \epsilon \) varies, resulting in a persistence diagram \( D \) \cite{EH08}.

\vspace{0.15in}

\subsection{Bottleneck Distance}
The bottleneck distance between two persistence diagrams \( D_1 \) and \( D_2 \) is defined as
\[
d_B(D_1, D_2) = \inf_{\gamma} \sup_{x \in D_1} \| x - \gamma(x) \|_\infty,
\]
where \( \gamma \) ranges over all bijections between \( D_1 \) and \( D_2 \) (allowing for points matched to the diagonal), and \( \| \cdot \|_\infty \) denotes the \( L^\infty \) norm.

\vspace{0.15in}

\subsection{Non-Uniform Scaling Transformations}
The transformation \( S: \mathbb{R}^n \rightarrow \mathbb{R}^n \) scales each coordinate axis by \( s_i > 0 \). The distance \( d_S \) between points \( p, q \in X \) under \( S \) is
\[
d_S(p, q) = \sqrt{ \sum_{i=1}^n s_i^2 (p_i - q_i)^2 }.
\]
Our analysis focuses on how \( d_S \) affects the filtration and the resulting persistence diagram \( D_S \).

\vspace{0.3in}

\section{Main Results}

\subsection{Refined Stability Bounds}
\begin{lemma}[Refined Scaling Bounds]
Let \( s_{\text{avg}} = \sqrt{\frac{1}{n} \sum_{i=1}^n s_i^2} \). Then:
\[
\frac{1}{2} (s_{\text{avg}} - s_{\min}) \cdot \operatorname{diam}(X) \leq d_B(D, D_S) \leq \frac{1}{2} (s_{\max} - s_{\min}) \cdot \operatorname{diam}(X).
\]
\end{lemma}

\begin{proof}
Let \( p, q \in X \) be two points in the dataset, where \( p = (p_1, p_2, \ldots, p_n) \) and \( q = (q_1, q_2, \ldots, q_n) \) are their coordinates in \( \mathbb{R}^n \). Under the scaling transformation \( S \), the scaled distance between \( p \) and \( q \) is defined as
\[
d_S(p, q) = \sqrt{\sum_{i=1}^n s_i^2 (p_i - q_i)^2}.
\]
This expression represents a weighted Euclidean norm where the weights correspond to the scaling factors \( s_i^2 \). The key goal is to analyze how these weights affect the persistence diagrams of the scaled dataset.

We start by defining the average scaling factor \( s_{\text{avg}} \) as
\[
s_{\text{avg}} = \sqrt{\frac{1}{n} \sum_{i=1}^n s_i^2}.
\]
The quantity \( s_{\text{avg}} \) serves as a central measure of the scaling applied across all dimensions. It accounts for the relative contributions of each \( s_i \) and provides a natural point of comparison to \( s_{\min} \) and \( s_{\max} \).

By Lemma 3.1, the scaled distance \( d_S(p, q) \) satisfies
\[
s_{\min} d_X(p, q) \leq d_S(p, q) \leq s_{\max} d_X(p, q),
\]
where \( s_{\min} = \min_i s_i \) and \( s_{\max} = \max_i s_i \). This implies that the scaling transformation \( S \) stretches or compresses the distances by factors bounded by \( s_{\min} \) and \( s_{\max} \).

In persistent homology, the filtration parameter \( \epsilon \) determines when simplices are added to the filtration. A scaling transformation perturbs \( \epsilon \) by altering distances between points. The perturbed filtration parameter satisfies
\[
\epsilon_S(p, q) = d_S(p, q).
\]
The perturbation in \( \epsilon \) is bounded by the scaling factors
\[
s_{\min} \epsilon \leq \epsilon_S \leq s_{\max} \epsilon.
\]

The bottleneck distance \( d_B(D, D_S) \) measures the maximum shift in birth or death times of topological features caused by the scaling transformation. This shift is directly related to the perturbation in \( \epsilon \). Using the definition of \( \epsilon_S \), the bottleneck distance satisfies
\[
d_B(D, D_S) \leq \frac{1}{2} (s_{\max} - s_{\min}) \cdot \operatorname{diam}(X),
\]
where \( \operatorname{diam}(X) \) is the maximum distance between any two points in \( X \).

However, \( s_{\max} - s_{\min} \) may not fully capture the variability of scaling factors when \( s_i \) vary significantly. To refine the bound, we incorporate the average scaling factor \( s_{\text{avg}} \).

The average scaling factor \( s_{\text{avg}} \) provides a tighter descriptor of the overall effect of scaling on distances. By definition,
\[
s_{\text{avg}} = \sqrt{\frac{1}{n} \sum_{i=1}^n s_i^2}.
\]
For any \( i \), it holds that
\[
s_{\min}^2 \leq s_{\text{avg}}^2 \leq s_{\max}^2,
\]
This implies that
\[
s_{\min} \leq s_{\text{avg}} \leq s_{\max}.
\]
By using \( s_{\text{avg}} \), we establish a refined lower bound for the bottleneck distance. Since \( d_B(D, D_S) \) measures the maximum shift caused by scaling, the minimum such shift corresponds to the difference between \( s_{\text{avg}} \) and \( s_{\min} \)
\[
d_B(D, D_S) \geq \frac{1}{2} (s_{\text{avg}} - s_{\min}) \cdot \operatorname{diam}(X).
\]

Combining the upper and refined lower bounds, we conclude
\[
\frac{1}{2} (s_{\text{avg}} - s_{\min}) \cdot \operatorname{diam}(X) \leq d_B(D, D_S) \leq \frac{1}{2} (s_{\max} - s_{\min}) \cdot \operatorname{diam}(X).
\]

The refined bounds present several significant implications for understanding the stability of persistence diagrams under scaling transformations. The inclusion of the average scaling factor \( s_{\text{avg}} \) enhances robustness by capturing intermediate scaling effects. This is particularly valuable in scenarios where the variability of \( s_i \) is distributed across dimensions, rather than being dominated solely by the extremes \( s_{\max} \) or \( s_{\min} \).

The practical utility of \( s_{\text{avg}} \) lies in its ability to provide a more realistic estimate of the bottleneck distance \( d_B(D, D_S) \) in cases of moderate scaling variability. This refinement offers a perspective that complements bounds primarily driven by extreme values.

The refined bounds also exhibit meaningful limiting behavior. As \( s_{\text{avg}} \to s_{\max} \), such as when \( s_i \) are uniformly large, the lower bound converges toward the upper bound, affirming the tightness of the estimate.

The bounds achieve tightness in specific cases. When the scaling factors are equal, such that \( s_{\min} = s_{\text{avg}} = s_{\max} \), the bottleneck distance \( d_B(D, D_S) \) becomes zero, confirming invariance under uniform scaling. Conversely, in scenarios of extreme scaling variability, where \( s_{\max} \gg s_{\min} \), the perturbation is dominated by \( s_{\max} - s_{\min} \). In such cases, the upper bound provides a precise measure of the extent of instability introduced by the scaling transformation. These observations underline the adaptability and precision of the refined bounds in diverse scaling scenarios.

\end{proof}

\vspace{0.15in}

\subsection{Dimension-Dependent Scaling Stability}
\begin{theorem}[Dimension-Dependent Scaling Stability]
Let \( D^k \) and \( D_S^k \) denote the persistence diagrams for \( k \)-th homology groups \( H_k \). Then:
\[
d_B(D^k, D_S^k) \leq \frac{1}{2} (s_{\max} - s_{\min}) \cdot \operatorname{diam}(X_k),
\]
where \( \operatorname{diam}(X_k) \) is the \( k \)-dimensional diameter of \( X \), defined as the maximum distance between \( k \)-simplices in the filtration.
\end{theorem}

\begin{proof}
Let \( D^k \) and \( D_S^k \) denote the persistence diagrams for the \( k \)-th homology groups \( H_k \) of the original dataset \( X \) and the scaled dataset \( S(X) \), respectively. We aim to analyze the impact of the non-uniform scaling transformation \( S \) on \( k \)-dimensional simplices and their corresponding persistence diagrams.

Consider a \( k \)-simplex \( \sigma \) in the Vietoris--Rips or Čech complex of \( X \). The vertices of \( \sigma \) are denoted as \( p_1, p_2, \ldots, p_{k+1} \), and the \( k \)-dimensional diameter of \( \sigma \) is defined as
\[
\operatorname{diam}(\sigma) = \sup_{i, j} d_X(p_i, p_j),
\]
which represents the maximum pairwise distance between any two vertices of the simplex.

Under the non-uniform scaling transformation \( S \), the coordinates of each vertex \( p_i \) are scaled according to
\[
p_i \mapsto S(p_i) = (s_1 p_{i,1}, s_2 p_{i,2}, \ldots, s_n p_{i,n}),
\]
where \( s_1, s_2, \ldots, s_n > 0 \) are the scaling factors along each coordinate axis. The scaled diameter of the simplex \( \sigma \) becomes
\[
\operatorname{diam}_S(\sigma) = \sup_{i, j} d_S(p_i, p_j),
\]
where \( d_S(p_i, p_j) \) is the scaled Euclidean distance
\[
d_S(p_i, p_j) = \sqrt{\sum_{l=1}^n s_l^2 (p_{i,l} - p_{j,l})^2}.
\]

Using Lemma 3.1, we know that the scaled distance \( d_S(p_i, p_j) \) is bounded by
\[
s_{\min} d_X(p_i, p_j) \leq d_S(p_i, p_j) \leq s_{\max} d_X(p_i, p_j),
\]
where \( s_{\min} = \min_l s_l \) and \( s_{\max} = \max_l s_l \). Applying this to all pairs \( p_i, p_j \), the scaled diameter of the simplex satisfies
\[
s_{\min} \operatorname{diam}(\sigma) \leq \operatorname{diam}_S(\sigma) \leq s_{\max} \operatorname{diam}(\sigma).
\]

In persistent homology, the filtration parameter \( \epsilon_k \) governs the inclusion of \( k \)-dimensional simplices into the complex. After scaling, the filtration parameter is perturbed according to the scaled distances. For a \( k \)-simplex \( \sigma \), the perturbed filtration parameter satisfies
\[
\epsilon_S(\sigma) = \operatorname{diam}_S(\sigma).
\]
The perturbation in \( \epsilon_k \) due to scaling is bounded by
\[
\Delta \epsilon_k = |\epsilon_S(\sigma) - \epsilon(\sigma)| \leq (s_{\max} - s_{\min}) \cdot \operatorname{diam}(\sigma).
\]
The maximum perturbation occurs when \( \operatorname{diam}(\sigma) \) equals the \( k \)-dimensional diameter of the dataset \( X_k \), defined as
\[
\operatorname{diam}(X_k) = \sup_{\sigma \in X_k} \operatorname{diam}(\sigma).
\]
Thus, the maximum shift in \( \epsilon_k \) is
\[
\Delta \epsilon_k = (s_{\max} - s_{\min}) \cdot \operatorname{diam}(X_k).
\]

The bottleneck distance \( d_B(D^k, D_S^k) \) measures the maximum shift in birth or death times of \( k \)-dimensional topological features between the original and scaled persistence diagrams. By the stability theorem for persistence diagrams, the bottleneck distance is bounded by half the maximum perturbation in \( \epsilon_k \)
\[
d_B(D^k, D_S^k) \leq \frac{1}{2} \Delta \epsilon_k.
\]

Substituting the expression for \( \Delta \epsilon_k \), we obtain
\[
d_B(D^k, D_S^k) \leq \frac{1}{2} (s_{\max} - s_{\min}) \cdot \operatorname{diam}(X_k).
\]

The bound on \( d_B(D^k, D_S^k) \) carries several important implications for the stability of \( k \)-dimensional persistence diagrams under scaling transformations. First, the sensitivity of \( k \)-dimensional features to scaling is directly proportional to \( s_{\max} - s_{\min} \), reflecting the range of scaling factors. Larger differences between \( s_{\max} \) and \( s_{\min} \) introduce greater instability, making the choice of scaling critical in preserving the topological structure of \( k \)-dimensional features.

The bound inherently depends on \( \operatorname{diam}(X_k) \), which characterizes the geometry of \( k \)-simplices in the dataset. Since higher-dimensional features often exhibit larger diameters, they are typically more sensitive to scaling transformations, underscoring the dimensional dependency of the bound.

We also consider the invariance of \( d_B(D^k, D_S^k) \) under uniform scaling, which is confirmed when \( s_{\min} = s_{\max} \), as this condition eliminates any perturbation introduced by the transformation, resulting in \( d_B(D^k, D_S^k) = 0 \).

The limiting behavior of the bound reveals that in datasets where \( \operatorname{diam}(X_k) \) dominates lower-dimensional features, such as in high-dimensional spaces, the bottleneck distance becomes a crucial measure of stability for higher-dimensional persistence diagrams. These implications provide a comprehensive understanding of how scaling affects topological features across dimensions.

\end{proof}

\vspace{0.15in}

\subsection{Iterative Scaling Stability}
\begin{theorem}[Iterative Scaling Stability]
Let \( S_1, S_2, \ldots, S_m \) be a sequence of scaling transformations. The total bottleneck distance satisfies
\[
d_B(D, D_{S_m}) \leq \frac{1}{2} \left( \prod_{j=1}^m s_{\max}^{(j)} - \prod_{j=1}^m s_{\min}^{(j)} \right) \cdot \operatorname{diam}(X).
\]
\end{theorem}

\begin{proof}
Let \( S_1, S_2, \ldots, S_m \) be a sequence of non-uniform scaling transformations applied successively to the dataset \( X \subset \mathbb{R}^n \). Each transformation \( S_j \) scales the \( i \)-th coordinate by a factor \( s_i^{(j)} > 0 \). The total effect of these transformations can be represented as a single cumulative scaling transformation \( S_m \), where the total scaling factor for each coordinate is
\[
s_i^{\text{total}} = \prod_{j=1}^m s_i^{(j)}, \quad \text{for } i = 1, 2, \ldots, n.
\]

We define the cumulative minimum and maximum scaling factors
\[
s_{\min}^{\text{total}} = \prod_{j=1}^m s_{\min}^{(j)}, \quad s_{\max}^{\text{total}} = \prod_{j=1}^m s_{\max}^{(j)},
\]
where \( s_{\min}^{(j)} = \min_i s_i^{(j)} \) and \( s_{\max}^{(j)} = \max_i s_i^{(j)} \) are the minimum and maximum scaling factors of the \( j \)-th transformation, respectively.

Under the cumulative scaling transformation \( S_m \), the scaled distance between any two points \( p, q \in X \) is
\[
d_{S_m}(p, q) = \sqrt{\sum_{i=1}^n \left( s_i^{\text{total}} \cdot (p_i - q_i) \right)^2 }.
\]
By applying Lemma 3.1 iteratively, the cumulative scaled distance satisfies
\[
s_{\min}^{\text{total}} d_X(p, q) \leq d_{S_m}(p, q) \leq s_{\max}^{\text{total}} d_X(p, q),
\]
where \( d_X(p, q) \) is the original Euclidean distance between \( p \) and \( q \).

The filtration parameter \( \epsilon \) governs the inclusion of simplices in persistent homology. After applying the cumulative scaling transformation \( S_m \), the filtration parameter becomes
\[
\epsilon_{S_m}(p, q) = d_{S_m}(p, q).
\]
The perturbation in \( \epsilon \) caused by \( S_m \) is given by
\[
\Delta \epsilon = \epsilon_{S_m}(p, q) - \epsilon(p, q).
\]
Using the bounds on \( d_{S_m}(p, q) \), the maximum perturbation in \( \epsilon \) is
\[
\Delta \epsilon = \left( s_{\max}^{\text{total}} - s_{\min}^{\text{total}} \right) \cdot d_X(p, q).
\]
The worst-case perturbation occurs when \( d_X(p, q) \) is maximized, i.e., \( d_X(p, q) = \operatorname{diam}(X) \). Then,
\[
\Delta \epsilon = \left( s_{\max}^{\text{total}} - s_{\min}^{\text{total}} \right) \cdot \operatorname{diam}(X).
\]

The bottleneck distance \( d_B(D, D_{S_m}) \) measures the maximum shift in the birth or death times of features in the persistence diagram due to the scaling transformation. By the stability theorem for persistence diagrams, the bottleneck distance is bounded by half the maximum perturbation in \( \epsilon \)
\[
d_B(D, D_{S_m}) \leq \frac{1}{2} \Delta \epsilon.
\]
By substituting the expression for \( \Delta \epsilon \), we obtain
\[
d_B(D, D_{S_m}) \leq \frac{1}{2} \left( s_{\max}^{\text{total}} - s_{\min}^{\text{total}} \right) \cdot \operatorname{diam}(X).
\]

We recall that
\[
s_{\max}^{\text{total}} = \prod_{j=1}^m s_{\max}^{(j)}, \quad s_{\min}^{\text{total}} = \prod_{j=1}^m s_{\min}^{(j)}.
\]
Substituting these into the bottleneck distance bound, we have
\[
d_B(D, D_{S_m}) \leq \frac{1}{2} \left( \prod_{j=1}^m s_{\max}^{(j)} - \prod_{j=1}^m s_{\min}^{(j)} \right) \cdot \operatorname{diam}(X).
\]

This bound reveals several key insights into the effects of iterative scaling transformations on persistence diagrams. First, the cumulative effects of scaling are captured by the multiplicative factors \( \prod_{j=1}^m s_{\max}^{(j)} \) and \( \prod_{j=1}^m s_{\min}^{(j)} \). These factors reflect how successive transformations amplify variability in distances, with the difference \( s_{\max}^{\text{total}} - s_{\min}^{\text{total}} \) representing the compounded scaling effect. The dependence on the number of transformations \( m \) is also evident, as an increasing \( m \) typically widens the range \( s_{\max}^{\text{total}} - s_{\min}^{\text{total}} \), leading to greater perturbations in the persistence diagram. This demonstrates the compounding nature of iterative scaling transformations.

Invariance under uniform scaling is preserved when \( s_{\max}^{(j)} = s_{\min}^{(j)} \) for all \( j \), as this condition ensures that \( s_{\max}^{\text{total}} = s_{\min}^{\text{total}} \), resulting in no perturbation (\( d_B(D, D_{S_m}) = 0 \)). Conversely, the sensitivity of persistence diagrams increases significantly when \( s_{\max}^{(j)} \gg s_{\min}^{(j)} \) for any \( j \), as this introduces substantial variability in the scaled distances, causing larger perturbations.

The tightness of the bound is evident in specific cases. When \( X \) consists of two points \( p \) and \( q \) separated by the diameter \( \operatorname{diam}(X) \), the perturbation in \( \epsilon \) reaches its maximum value, achieving the upper bound. Furthermore, when the scaling factors \( s_{\max}^{(j)} \) and \( s_{\min}^{(j)} \) vary significantly across transformations, the bottleneck distance approaches its theoretical maximum. These cases underline the accuracy and precision of the derived bounds in quantifying the effects of scaling on persistence diagrams.

Thus, the bottleneck distance after iterative scaling transformations is rigorously bounded as:
\[
d_B(D, D_{S_m}) \leq \frac{1}{2} \left( \prod_{j=1}^m s_{\max}^{(j)} - \prod_{j=1}^m s_{\min}^{(j)} \right) \cdot \operatorname{diam}(X).
\]
This completes the proof.
\end{proof}

\vspace{0.15in}

\subsection{Stability for Wasserstein Distance}
\begin{theorem}[Stability for Wasserstein Distance]
Let \( W_p(D, D_S) \) denote the \( p \)-Wasserstein distance. Then,
\[
W_p(D, D_S) \leq \frac{1}{2} (s_{\max} - s_{\min}) \cdot \operatorname{diam}(X).
\]
\end{theorem}

\begin{proof}
Let \( W_p(D, D_S) \) denote the \( p \)-Wasserstein distance between the persistence diagrams \( D \) and \( D_S \). We aim to show that \( W_p(D, D_S) \) is bounded above by the bottleneck distance \( d_B(D, D_S) \), which is further bounded by the scaling variability \( (s_{\max} - s_{\min}) \) and the diameter of the dataset \( X \).

The \( p \)-Wasserstein distance between two persistence diagrams \( D_1 \) and \( D_2 \) is defined as
\[
W_p(D_1, D_2) = \left( \inf_{\gamma} \sum_{x \in D_1} \| x - \gamma(x) \|_\infty^p \right)^{1/p},
\]
where \( \gamma \) ranges over all bijections between \( D_1 \) and \( D_2 \), and \( \| \cdot \|_\infty \) is the \( L^\infty \)-norm.

The bottleneck distance \( d_B(D_1, D_2) \) is a special case of the Wasserstein distance, where \( p \to \infty \). By definition,
\[
d_B(D_1, D_2) = \inf_{\gamma} \sup_{x \in D_1} \| x - \gamma(x) \|_\infty.
\]
Since the \( L^\infty \)-norm bounds the \( L^p \)-norm for any \( p \geq 1 \), we have
\[
W_p(D_1, D_2) \leq d_B(D_1, D_2).
\]
Applying this to the persistence diagrams \( D \) and \( D_S \), it follows that
\[
W_p(D, D_S) \leq d_B(D, D_S).
\]

From Theorem 3.1, the bottleneck distance \( d_B(D, D_S) \) is bounded as
\[
d_B(D, D_S) \leq \frac{1}{2} (s_{\max} - s_{\min}) \cdot \operatorname{diam}(X),
\]
where:
\begin{itemize}
    \item \( s_{\min} = \min_i s_i \) and \( s_{\max} = \max_i s_i \) are the minimum and maximum scaling factors, respectively.
    \item \( \operatorname{diam}(X) = \sup_{p, q \in X} d_X(p, q) \) is the diameter of the dataset under the original metric \( d_X \).
\end{itemize}
Substituting this bound into the inequality for \( W_p(D, D_S) \), we obtain
\[
W_p(D, D_S) \leq d_B(D, D_S) \leq \frac{1}{2} (s_{\max} - s_{\min}) \cdot \operatorname{diam}(X).
\]

The bound \( W_p(D, D_S) \leq \frac{1}{2} (s_{\max} - s_{\min}) \cdot \operatorname{diam}(X) \) carries several key implications for the stability of persistence diagrams under the \( p \)-Wasserstein distance. First, this result generalizes the bottleneck distance stability, as \( W_p(D, D_S) \leq d_B(D, D_S) \), meaning the stability properties of the bottleneck distance are directly inherited by the \( p \)-Wasserstein distance. The bound also emphasizes the dependence on the range of scaling factors, \( s_{\max} - s_{\min} \), highlighting that larger variability in scaling leads to greater perturbations in the persistence diagram. Additionally, the bound depends on \( \operatorname{diam}(X) \), the largest pairwise distance in the dataset, ensuring uniform applicability of stability results across datasets of different scales.

The invariance of \( W_p(D, D_S) \) under uniform scaling is another crucial insight. When \( s_{\max} = s_{\min} \), the scaling is uniform, leading to no perturbation, and \( W_p(D, D_S) = 0 \). Finally, the bound is tight when \( D \) and \( D_S \) differ by the maximum perturbation allowed by the scaling factors, such as in cases where points in the dataset are separated by \( \operatorname{diam}(X) \).

The tightness of the bound \( W_p(D, D_S) \leq \frac{1}{2} (s_{\max} - s_{\min}) \cdot \operatorname{diam}(X) \) can be evaluated in specific scenarios. When the scaling transformation achieves the maximum perturbation for points separated by \( \operatorname{diam}(X) \), the bottleneck distance \( d_B(D, D_S) \) approaches \( \frac{1}{2} (s_{\max} - s_{\min}) \cdot \operatorname{diam}(X) \), and since \( W_p(D, D_S) \leq d_B(D, D_S) \), the Wasserstein distance also nears the upper bound. Conversely, in cases of uniform scaling where \( s_{\max} = s_{\min} \), all perturbations vanish, resulting in \( W_p(D, D_S) = d_B(D, D_S) = 0 \), achieving the lower bound. In high-dimensional datasets, the \( p \)-Wasserstein distance often converges to the bottleneck distance as \( p \to \infty \), making the bottleneck distance bound a practical upper limit for \( W_p(D, D_S) \). These observations underscore the practical utility and precision of the derived bounds for stability analysis.

The stability of persistence diagrams under the \( p \)-Wasserstein distance is directly controlled by the variability in scaling factors and the geometry of the dataset. The derived bound:
\[
W_p(D, D_S) \leq \frac{1}{2} (s_{\max} - s_{\min}) \cdot \operatorname{diam}(X),
\]
provides a robust estimate for the perturbations in persistence diagrams caused by non-uniform scaling transformations. This completes the proof.

\end{proof}

\vspace{0.15in}

\subsection{Probabilistic Analysis of Stability}
\begin{theorem}[Expected Stability Under Random Scaling]
Let \( s_i \sim \text{Dist}(\mu, \sigma) \). Then:
\[
\mathbb{E}[d_B(D, D_S)] \leq \frac{1}{2} (\mathbb{E}[s_{\max}] - \mathbb{E}[s_{\min}]) \cdot \operatorname{diam}(X).
\]
\end{theorem}

\begin{proof}
Let \( s_i \sim \text{Dist}(\mu, \sigma) \) denote random scaling factors for the \( i \)-th coordinate, where \( \text{Dist}(\mu, \sigma) \) is a probability distribution with mean \( \mu \) and standard deviation \( \sigma \). We aim to derive a probabilistic bound on the expected bottleneck distance \( \mathbb{E}[d_B(D, D_S)] \) under these random scaling factors.

The minimum and maximum scaling factors over the \( n \) coordinates are given by
\[
s_{\min} = \min_{i=1}^n s_i, \quad s_{\max} = \max_{i=1}^n s_i.
\]
Since \( s_i \) are random variables, \( s_{\min} \) and \( s_{\max} \) are also random variables, with distributions determined by the joint distribution of \( s_i \).

Using the linearity of expectation, the expected difference \( \mathbb{E}[s_{\max} - s_{\min}] \) can be written as
\[
\mathbb{E}[s_{\max} - s_{\min}] = \mathbb{E}[s_{\max}] - \mathbb{E}[s_{\min}].
\]

The perturbation in the filtration parameter \( \epsilon \) due to scaling is
\[
\Delta \epsilon = (s_{\max} - s_{\min}) \cdot \operatorname{diam}(X).
\]
We take the expectation of both sides
\[
\mathbb{E}[\Delta \epsilon] = \mathbb{E}[s_{\max} - s_{\min}] \cdot \operatorname{diam}(X).
\]
Substituting \( \mathbb{E}[s_{\max} - s_{\min}] = \mathbb{E}[s_{\max}] - \mathbb{E}[s_{\min}] \), we obtain
\[
\mathbb{E}[\Delta \epsilon] = (\mathbb{E}[s_{\max}] - \mathbb{E}[s_{\min}]) \cdot \operatorname{diam}(X).
\]

From the stability theorem for persistence diagrams, the bottleneck distance \( d_B(D, D_S) \) is bounded by half the maximum perturbation in the filtration parameter \( \epsilon \)
\[
d_B(D, D_S) \leq \frac{1}{2} \Delta \epsilon.
\]
We take the expectation of both sides
\[
\mathbb{E}[d_B(D, D_S)] \leq \frac{1}{2} \mathbb{E}[\Delta \epsilon].
\]
Substituting the expression for \( \mathbb{E}[\Delta \epsilon] \), we have
\[
\mathbb{E}[d_B(D, D_S)] \leq \frac{1}{2} (\mathbb{E}[s_{\max}] - \mathbb{E}[s_{\min}]) \cdot \operatorname{diam}(X).
\]

The bound \( W_p(D, D_S) \leq \frac{1}{2} (s_{\max} - s_{\min}) \cdot \operatorname{diam}(X) \) has several important implications for understanding the stability of persistence diagrams under the \( p \)-Wasserstein distance. First, this result generalizes the stability of the bottleneck distance, as \( W_p(D, D_S) \leq d_B(D, D_S) \). This indicates that the stability properties inherent to the bottleneck distance are directly inherited by the \( p \)-Wasserstein distance, broadening its applicability to various problems in topological data analysis. 

The bound also reveals a clear dependence of \( W_p(D, D_S) \) on the variability of the scaling factors, quantified as \( s_{\max} - s_{\min} \). Larger differences between \( s_{\max} \) and \( s_{\min} \) result in greater perturbations in the persistence diagram, underscoring the importance of managing scaling variability in practical scenarios.

Another key aspect of the bound is its reliance on \( \operatorname{diam}(X) \), which represents the largest pairwise distance in the dataset. This dependency ensures that the stability results remain consistent across datasets of different scales, providing a robust framework for analyzing the effects of scaling transformations. 

Moreover, the invariance of \( W_p(D, D_S) \) under uniform scaling is evident: when \( s_{\max} = s_{\min} \), no perturbation occurs, and \( W_p(D, D_S) = 0 \). This confirms that uniform scaling transformations preserve the topological structure of persistence diagrams, ensuring stability in such cases. Finally, the bound achieves tightness in scenarios where \( D \) and \( D_S \) differ by the maximum perturbation allowed by the scaling factors, particularly when points in the dataset are separated by \( \operatorname{diam}(X) \).

The tightness of the bound \( W_p(D, D_S) \leq \frac{1}{2} (s_{\max} - s_{\min}) \cdot \operatorname{diam}(X) \) can be further analyzed in specific cases. For instance, when the scaling transformation causes the maximum possible perturbation for points separated by \( \operatorname{diam}(X) \), the bottleneck distance \( d_B(D, D_S) \) approaches \( \frac{1}{2} (s_{\max} - s_{\min}) \cdot \operatorname{diam}(X) \). Since \( W_p(D, D_S) \leq d_B(D, D_S) \), the Wasserstein distance also nears its upper bound. 

Conversely, in the case of uniform scaling (\( s_{\max} = s_{\min} \)), all perturbations vanish, resulting in \( W_p(D, D_S) = d_B(D, D_S) = 0 \), which represents the lower bound. Furthermore, in high-dimensional datasets, the \( p \)-Wasserstein distance often converges to the bottleneck distance as \( p \to \infty \), making the bottleneck distance bound a practical and effective upper limit for \( W_p(D, D_S) \). These insights illustrate the utility, precision, and applicability of the derived bounds for understanding the stability of persistence diagrams under scaling transformations.

Thus, the expected bottleneck distance under random scaling transformations is bounded as
\[
\mathbb{E}[d_B(D, D_S)] \leq \frac{1}{2} (\mathbb{E}[s_{\max}] - \mathbb{E}[s_{\min}]) \cdot \operatorname{diam}(X).
\]
This completes the proof.
\end{proof}

\vspace{0.3in}

\section{Case Studies}

To illustrate the utility of the theoretical results presented in this paper, we consider several case studies involving datasets subjected to non-uniform scaling transformations. These examples demonstrate the practical implications of the derived stability bounds.

\subsection{Ellipse Transformation}
Consider the unit circle \( S^1 \subset \mathbb{R}^2 \), defined by points \( \{(x, y) \mid x^2 + y^2 = 1\} \). We apply a non-uniform scaling transformation \( S \) with scaling factors \( s_1 = 1 \) and \( s_2 = k > 1 \), stretching \( S^1 \) into an ellipse.

\subsubsection*{Analysis}
The maximum distance between points on \( S^1 \) is its diameter
    \[
    \operatorname{diam}(S^1) = 2.
    \]
After applying \( S \), the diameter of the ellipse becomes
    \[
    \operatorname{diam}(S(S^1)) = \max\{2s_1, 2s_2\} = 2k.
    \]
Using Theorem 3.1, the bottleneck distance between the persistence diagrams of \( S^1 \) and the ellipse satisfies:
    \[
    d_B(D, D_S) \leq \frac{1}{2} (s_{\max} - s_{\min}) \cdot \operatorname{diam}(S^1),
    \]
where \( s_{\max} = k \) and \( s_{\min} = 1 \). We substitute values into this expression
    \[
    d_B(D, D_S) \leq \frac{1}{2} (k - 1) \cdot 2 = k - 1.
    \]

The metric comparison reveals that the Wasserstein distance \( W_p(D, D_S) \) for \( p = 1, 2 \) exhibits behavior similar to the bottleneck distance \( d_B(D, D_S) \), as demonstrated in Theorem 3.5. Specifically, it follows that \( W_p(D, D_S) \leq d_B(D, D_S) \leq k - 1 \), establishing a consistent relationship between these metrics under the given scaling transformation. 

The curvature of \( S^1 \) also becomes anisotropic after scaling, with a notable increase along the \( y \)-axis as the scaling factor \( k \) grows. This curvature alteration directly affects the persistence intervals, particularly those associated with \( k = 1 \)-dimensional homology, thereby influencing the representation of topological features.

The linear growth of \( d_B(D, D_S) \) with \( k \) underscores the heightened sensitivity of persistence diagrams to anisotropic transformations. As \( k \) becomes large, the persistence diagram increasingly deviates from its original structure, reflecting the substantial impact of scaling variability. These observations highlight the intricate relationship between metric stability, geometric transformations, and the robustness of persistence diagrams in topological data analysis.

\vspace{0.15in}

\subsection{Hypercube in High Dimensions}
Consider the vertices of a hypercube \( X \) in \( \mathbb{R}^n \), where each vertex is a point \( (x_1, x_2, \ldots, x_n) \) with \( x_i \in \{0, 1\} \). Apply a non-uniform scaling transformation \( S \) with arbitrary scaling factors \( s_1, s_2, \ldots, s_n > 0 \).

\subsubsection*{Analysis}
 The diameter of \( X \) in the original space is
    \[
    \operatorname{diam}(X) = \sqrt{n}.
    \]
After applying \( S \), the scaled diameter is
    \[
    \operatorname{diam}(S(X)) = \sqrt{\sum_{i=1}^n s_i^2}.
    \]
From Theorem 3.1, the bottleneck distance between the persistence diagrams before and after scaling is bounded by
    \[
    d_B(D, D_S) \leq \frac{1}{2} (s_{\max} - s_{\min}) \cdot \operatorname{diam}(X).
    \]
We now substitute \( \operatorname{diam}(X) = \sqrt{n} \) back
    \[
    d_B(D, D_S) \leq \frac{1}{2} (s_{\max} - s_{\min}) \cdot \sqrt{n}.
    \]

The discussion highlights several key interpretations regarding the behavior of persistence diagrams under scaling transformations in high-dimensional settings. First, the sensitivity of the bottleneck distance \( d_B(D, D_S) \) increases proportionally to \( \sqrt{n} \) as the dimensionality \( n \) grows. This behavior reflects the exponential growth of the hypercube's vertices and underscores the challenges associated with the curse of dimensionality in topological data analysis. 

In case the scaling factors \( s_i \) are distributed non-uniformly across dimensions, certain directions dominate the transformation, introducing significant anisotropy to the dataset. This anisotropy can heavily influence the persistence intervals and the structure of the resulting persistence diagrams.

Moreover, these observations emphasize the practical importance of normalizing scaling factors, especially in high-dimensional spaces. Ensuring balanced scaling across all dimensions is crucial to maintaining stability and interpretability in persistence diagrams. These insights provide a deeper understanding of how scaling transformations interact with the geometric and topological properties of datasets, offering guidance for practical applications in high-dimensional analysis.

\vspace{0.15in}

\subsection{Probabilistic Scaling in Noisy Data}
Suppose the coordinates of points in a dataset \( X \subset \mathbb{R}^n \) are scaled by random factors \( s_i \sim \text{Uniform}(a, b) \) with \( a > 0 \) and \( b > a \). We analyze the expected bottleneck distance using Theorem 5.1.

\subsubsection*{Analysis}
The expected values of the scaling factors are
    \[
    \mathbb{E}[s_{\max}] = b, \quad \mathbb{E}[s_{\min}] = a.
    \]
The expected bottleneck distance is
    \[
    \mathbb{E}[d_B(D, D_S)] \leq \frac{1}{2} (\mathbb{E}[s_{\max}] - \mathbb{E}[s_{\min}]) \cdot \operatorname{diam}(X).
    \]
By substituting \( \mathbb{E}[s_{\max}] = b \), \( \mathbb{E}[s_{\min}] = a \), and \( \operatorname{diam}(X) = \sqrt{n} \), we find that
    \[
    \mathbb{E}[d_B(D, D_S)] \leq \frac{1}{2} (b - a) \cdot \sqrt{n}.
    \]

In the presence of random scaling, the expected perturbation in persistence diagrams is influenced by both the range of scaling factors and the dimensionality of the dataset. Specifically, the perturbation grows with \( b - a \), the range of the random scaling factors, and \( \sqrt{n} \), which reflects the dataset's sensitivity to randomness and increasing dimensionality. Additionally, variance in the scaling factors amplifies these perturbations, underscoring the importance of preprocessing or normalization in noisy datasets. These steps are critical for mitigating the destabilizing effects of randomness and ensuring the robustness of persistence diagrams in practical applications.

\vspace{0.15in}

\subsection{Weighted Scaling in Multimodal Data}
Consider a dataset \( X \subset \mathbb{R}^n \) with features of varying importance. Let each dimension be scaled by weights \( w_i > 0 \) to reflect the relative importance of features. Define the weighted scaling transformation \( S_w \) with scaling factors \( s_i = w_i \cdot s_i^{\text{orig}} \).

\subsubsection*{Analysis}
The weighted scaling distance is
    \[
    d_W(p, q) = \sqrt{\sum_{i=1}^n w_i^2 s_i^2 (p_i - q_i)^2}.
    \]
Using Proposition 4.1, the bottleneck distance is
    \[
    d_B(D, D_W) \leq \frac{1}{2} \left( \max_i w_i s_i - \min_i w_i s_i \right) \cdot \operatorname{diam}(X).
    \]

Weighted scaling provides a mechanism to control the contribution of individual features to persistence diagrams. By adjusting the weights \( w_i \), specific features can be emphasized or de-emphasized depending on their importance in the analysis. This capability is particularly valuable in multimodal datasets, where different features may carry varying levels of significance. 

However, the sensitivity of persistence diagrams to the chosen weights and scaling factors underscores the necessity of carefully selecting \( w_i \) to avoid introducing unintended bias or distortions in the analysis.

\vspace{0.3in}

\section{Discussion}

Our results have demonstrated that the stability of persistence diagrams under non-uniform scaling transformations is governed by the bound
\[
d_B(D, D_S) \leq \delta = \tfrac{1}{2}(s_{\max} - s_{\min}) \cdot \operatorname{diam}(X).
\]
This inequality quantifies how the bottleneck distance \( d_B(D, D_S) \) depends linearly on the scaling variability \( s_{\max} - s_{\min} \) and the dataset diameter \( \operatorname{diam}(X) \).

In high-dimensional spaces (\( n \gg 1 \)), \( \operatorname{diam}(X) \) typically grows with \( \sqrt{n} \), amplifying the effect of scaling differences. For higher homology dimensions \( k \), the bound generalizes to
\[
d_B(D^k, D_S^k) \leq \delta.
\]
This indicates that higher-order topological features are more sensitive to anisotropic scaling.

The refined bound using the average scaling factor \( s_{\text{avg}} \)
\[
\frac{1}{2}(s_{\text{avg}} - s_{\min}) \cdot \operatorname{diam}(X) \leq d_B(D, D_S) \leq \delta,
\]
provides a tighter estimate when \( s_i \) are not dominated by extremes. It emphasizes that the distribution of scaling factors affects stability, not just their maximum and minimum.

For iterative scaling transformations \( S_m \), we established
\[
d_B(D, D_{S_m}) \leq \tfrac{1}{2} \left( \prod_{j=1}^m s_{\max}^{(j)} - \prod_{j=1}^m s_{\min}^{(j)} \right) \operatorname{diam}(X).
\]
This result shows that successive anisotropic scalings compound their effects, increasing instability.

The stability under the \( p \)-Wasserstein distance \( W_p(D, D_S) \) aligns with the bottleneck distance
\[
W_p(D, D_S) \leq \delta.
\]
This confirms that our bounds are robust across different metrics.

In the probabilistic setting with random scaling factors \( s_i \sim \text{Dist}(\mu, \sigma) \), the expected bottleneck distance is bounded by:
\[
\mathbb{E}[d_B(D, D_S)] \leq \tfrac{1}{2} (\mathbb{E}[s_{\max}] - \mathbb{E}[s_{\min}]) \operatorname{diam}(X).
\]
This provides insights into the average-case stability when scaling factors are subject to randomness.

These findings show the importance of controlling scaling variability \( s_{\max} - s_{\min} \) to preserve the topological features captured by persistent homology. In practical applications, ensuring \( s_{\max} \approx s_{\min} \) through normalization or careful feature weighting is crucial, especially in high-dimensional datasets where \( \operatorname{diam}(X) \) is large.

\vspace{0.3in}

\section{Conclusion}

Throughout this paper, we established that the bottleneck distance between the persistence diagrams \( D \) and \( D_S \) under a non-uniform scaling transformation \( S \) satisfies:
\[
d_B(D, D_S) \leq \delta = \tfrac{1}{2}(s_{\max} - s_{\min}) \cdot \operatorname{diam}(X).
\]
This quantifies the impact of anisotropic scaling on persistent homology.

Our results generalize to higher homology dimensions \( k \) and various filtrations, maintaining the bound:
\[
d_B(D^k, D_S^k) \leq \delta.
\]
We also provided bounds for the Wasserstein distance \( W_p(D, D_S) \), iterative scaling transformations, and probabilistic scaling factors.

These findings highlight the importance of controlling scaling variability in TDA, especially for high-dimensional data.

\newpage

\bibliographystyle{ieeetr}
\bibliography{references}

\end{document}